\documentclass[12pt]{article}
\usepackage{mathtools} 
\usepackage{amsmath, amssymb} 
\usepackage{amsthm}
\usepackage{amsrefs}
\usepackage{url} % formats url's nicely using \url{}
\usepackage{enumitem} 

\usepackage{graphicx}
\usepackage{tikz}
\usepackage{tkz-graph}
\usetikzlibrary{shapes}
\usetikzlibrary{arrows}
\usetikzlibrary{decorations.markings}
\usepackage{float}
\usepackage{graphicx}
\usepackage{caption,subcaption}

\DeclarePairedDelimiter\abs{\lvert}{\rvert}%
\DeclarePairedDelimiter\norm{\lVert}{\rVert}%

\makeatletter
\let\oldabs\abs
\def\abs{\@ifstar{\oldabs}{\oldabs*}}
\let\oldnorm\norm
\def\norm{\@ifstar{\oldnorm}{\oldnorm*}}
\makeatother

\newtheoremstyle{plainsl}%
	{\topsep}
	{\topsep}
	{\slshape} % only non-default setting
	{}
	{\normalfont\bfseries}
	{.}
	{ }
	{}

\swapnumbers

{\theoremstyle{plainsl}
\newtheorem{theorem}{Theorem}[section]
\newtheorem{lemma}[theorem]{Lemma}
\newtheorem{corollary}[theorem]{Corollary}}
{\theoremstyle{remark}
}

\renewcommand\proof{\noindent\textsl{Proof. }}
\newcommand\sqr[2]{{\vbox{\hrule height.#2pt
    \hbox{\vrule width.#2pt height#1pt \kern#1pt
        \vrule width.#2pt}\hrule height.#2pt}}}
\renewcommand\qed{%
	\ifmmode\eqno\sqr53
	\else\nolinebreak\ \hfill\sqr53\medbreak\fi}

\DeclareMathOperator{\rk}{rk}
\DeclareMathOperator{\tr}{tr}
\DeclareMathOperator{\col}{Col}

\DeclareMathOperator\spec{Spec}
\newcommand{\widehatP}{\widehat{P}}
\newcommand\scalemath[2]{\scalebox{#1}{\mbox{\ensuremath{\displaystyle #2}}}}

\usepackage{blkarray}

\usepackage{hyperref}

\title{Periodicity of bipartite walks on biregular graphs with conditional spectra }
\author{Qiuting Chen\\
 Department of Combinatorics \& Optimization,\\ University of Waterloo, Waterloo, Ontario, Canada\\ q243chen@uwaterloo.ca}
\begin{document}
\maketitle

\abstract{ In this paper we study a class of discrete quantum walks, known as bipartite walks. These include the well-known Grover's walks. Any discrete quantum walk is given by the powers of a unitary matrix $U$ indexed by arcs or edges of the underlying graph.  The walk is periodic if $U^k=I$ for some positive integer $k$. Kubota has given a characterization of periodicity of Grover's walk when the walk is defined on a regular bipartite graph with at most five eigenvalues. We extend Kubota's results\textemdash if a biregular graph $G$ has eigenvalues whose squares are algebraic integers with degree at most two, we characterize periodicity of the bipartite walk over $G$ in terms of its spectrum.  We apply periodicity results of bipartite walks to get a characterization of periodicity of Grover's walk on regular graphs.}

\section{Introduction}
Quantum walks are a quantum mechanical analogue of classical random walks. They are powerful tools for improving existing algorithms and developing new algorithms~\cite{szegedy2004,Childs2002}. Many properties of quantum walks have been studied in the context of both continuous and discrete quantum walks~\cite{avm_cts,Godsil2010,HadamardPED,PEDMIXED}. 

In this paper, we study periodicity of a discrete quantum walk model, called bipartite walks~\cite{ham}, which can be described as follows. Here is a description of bipartite walks. Given a bipartite graph $G$ with two color classes $C_0,C_1$, there are two partitions $\pi_0,\pi_1$ of $E(G)$, the edges of $G$. If two edges have a common vertex in $C_0$, then they are in the same cell in partition $\pi_0$. Similarly, if two edges have a common vertex in $C_1$, then they are in the same cell in partition $\pi_1$. For partitions $\pi_0$ and $\pi_1$, we define two projections $P$ and $Q$ onto the functions on $E(G)$ that are constant on cells of $\pi_0, \pi_1$ respectively. The transition matrix of the bipartite walk defined over $G$ is 
\[
	U=(2P-I)(2Q-I),
\]
which is a product of reflection $2P-I$ and reflection $2Q-I$.

Bipartite walks include many of commonly used discrete quantum walk models and so, the results of bipartite walks can be applied to many known models. This is one of reasons why we choose to study bipartite walks. Vertex-face walks~\cite{qwonembeddings} and Grover's walks are among the special cases of bipartite walks.

Periodicity is one of the major problems in study of quantum walks. 
In continuous quantum walks, the subject of periodicity is discussed using different Hamiltonians and various forms of states~\cite{periodicgraphs,PairState,GroupState}. Periodicity of discrete quantum walks is also studied using many different discrete quantum walk models~\cite{pedDQWonfiniteG,Barr}.
Periodic quantum walks help to design new quantum algorithms in quantum cryptology~\cite{orderchaos}. It is also of interest in development of quantum chaos control theory~\cite{Whaley_2015}.

In this paper, we study the periodicity of bipartite walks. Let $U$ denote the transition matrix of a bipartite walk. When there is a positive integer $k$ such that 
\[
U^k=I,\]
we say the walk is periodic. 

We say a graph is a biregular graph if it is a bipartite graph and all the vertices of the same class has the same degree. We give a characterization of periodicity of bipartite walks on the biregular graphs that squares of their eigenvalues are algebraic integer with degree at most two. The characterization is in terms of the spectrum of the graph the bipartite walk is defined on.

Recently in~\cite{ShoPED}, Kubota studies the periodicity of Grover's walk on regular bipartite graphs with at most five distinct eigenvalues. Since the number of distinct eigenvalues of graphs are at most five, it implies all the eigenvalues of graphs of interest in~\cite{ShoPED} are algebraic integers with degree at most two. Later in the paper, we show that Grover's walk is a special case of bipartite walk. Given a graph $G$, the transition matrix of the Grover's walk defined over $G$ is the exactly same as the transition matrix of the bipartite walk defined over the subdivision graph of $G$. So we can apply the result we have on bipartite walks to Grover's walk. 

Results in this paper extend the periodicity results of Grover's walks proved by Kubota in~\cite{ShoPED}. More specifically, results in this paper remove the constraint on the number of distinct eigenvalues and we do not require the underlying regular graph to be bipartite. We only requires that the underlying graph is regular and it has eigenvalues whose squares are algebraic integers with degree at most two. There will be examples of periodic graphs shown in the paper that are not explained by the results in~\cite{ShoPED}.

Moreover, we show that given a regular bipartite graph $G$, if the bipartite walk defined on $G$ is periodic with period $k$, then the Grover's walk defined on $G$ is also periodic and it has period $2k$. Over certain graphs, bipartite walk model has a advantage over Grover's walk model, which hopefully gives people more motivation to study bipartite walks.

\section{Bipartite walks}
\label{BW}
Let $G$ be a bipartite graph with two color classes $C_0,C_1$. 
Now we define two partitions of the edges of $G$, denoted by $\pi_0,\pi_1$ respectively.  If two edges have the same end $x$ in $C_0$, then they belong to the same cell of $\pi_0$. Similarly, if two edges have the same end $y$ in $C_1$, then they belong to the same cell of $\pi_1$.

Given a matrix $M$, we \textsl{normalize} it by scaling each column of $M$ to a unit vector.
Let $P_0,P_1$ be characteristic matrix of $\pi_0,\pi_1$ respectively and let $\widehatP_0,\widehatP_1$ denote the normalized $P_0,P_1$ respectively.

Let 
\[
	P=\widehatP_0\widehatP_0^T,\quad Q=\widehatP_1\widehatP_1^T
\] 
be the projections onto the vectors that is constant on the cells of $\pi_0,\pi_1$ respectively. It is easy to check that $2P-I$ and $2Q-I$ are two reflections. We define the transition matrix of the bipartite walk over $G$ to be the product of two reflections, i.e.,
\[
	U=\left(2P-I\right)\left(2Q-I\right).
\] 

\begin{figure}[H]
\begin{center}

\begin{tikzpicture}
\definecolor{cv0}{rgb}{0.0,0.0,0.0}
\definecolor{cfv0}{rgb}{1.0,1.0,1.0}
\definecolor{clv0}{rgb}{0.0,0.0,0.0}
\definecolor{cv1}{rgb}{0.0,0.0,0.0}
\definecolor{cfv1}{rgb}{1.0,1.0,1.0}
\definecolor{clv1}{rgb}{0.0,0.0,0.0}
\definecolor{cv2}{rgb}{0.0,0.0,0.0}
\definecolor{cfv2}{rgb}{1.0,1.0,1.0}
\definecolor{clv2}{rgb}{0.0,0.0,0.0}
\definecolor{cv3}{rgb}{0.0,0.0,0.0}
\definecolor{cfv3}{rgb}{1.0,1.0,1.0}
\definecolor{clv3}{rgb}{0.0,0.0,0.0}
\definecolor{cv4}{rgb}{0.0,0.0,0.0}
\definecolor{cfv4}{rgb}{1.0,1.0,1.0}
\definecolor{clv4}{rgb}{0.0,0.0,0.0}
\definecolor{cv5}{rgb}{0.0,0.0,0.0}
\definecolor{cfv5}{rgb}{1.0,1.0,1.0}
\definecolor{clv5}{rgb}{0.0,0.0,0.0}
\definecolor{cv6}{rgb}{0.0,0.0,0.0}
\definecolor{cfv6}{rgb}{1.0,1.0,1.0}
\definecolor{clv6}{rgb}{0.0,0.0,0.0}
\definecolor{cv7}{rgb}{0.0,0.0,0.0}
\definecolor{cfv7}{rgb}{1.0,1.0,1.0}
\definecolor{clv7}{rgb}{0.0,0.0,0.0}
\definecolor{cv0v1}{rgb}{0.0,0.0,0.0}
\definecolor{cv0v5}{rgb}{0.0,0.0,0.0}
\definecolor{cv1v2}{rgb}{0.0,0.0,0.0}
\definecolor{cv1v4}{rgb}{0.0,0.0,0.0}
\definecolor{cv2v3}{rgb}{0.0,0.0,0.0}
\definecolor{cv5v6}{rgb}{0.0,0.0,0.0}
\definecolor{cv6v7}{rgb}{0.0,0.0,0.0}
\Vertex[style={minimum size=1.0cm,draw=cv0,fill=cfv0,text=clv0,shape=circle},LabelOut=false,L=\hbox{$0$},x=0cm,y=4cm]{v0}
\Vertex[style={minimum size=1.0cm,draw=cv1,fill=cfv1,text=clv1,shape=circle},LabelOut=false,L=\hbox{$1$},x=4cm,y=4cm]{v1}
\Vertex[style={minimum size=1.0cm,draw=cv2,fill=cfv2,text=clv2,shape=circle},LabelOut=false,L=\hbox{$2$},x=0cm,y=3cm]{v2}
\Vertex[style={minimum size=1.0cm,draw=cv3,fill=cfv3,text=clv3,shape=circle},LabelOut=false,L=\hbox{$3$},x=4cm,y=3cm]{v3}
\Vertex[style={minimum size=1.0cm,draw=cv4,fill=cfv4,text=clv4,shape=circle},LabelOut=false,L=\hbox{$4$},x=0.0cm,y=2cm]{v4}
\Vertex[style={minimum size=1.0cm,draw=cv5,fill=cfv5,text=clv5,shape=circle},LabelOut=false,L=\hbox{$5$},x=4cm,y=2cm]{v5}
\Vertex[style={minimum size=1.0cm,draw=cv6,fill=cfv6,text=clv6,shape=circle},LabelOut=false,L=\hbox{$6$},x=0cm,y=1cm]{v6}
\Vertex[style={minimum size=1.0cm,draw=cv7,fill=cfv7,text=clv7,shape=circle},LabelOut=false,L=\hbox{$7$},x=4cm,y=1.0cm]{v7}
\Edge[lw=0.1cm,style={color=cv0v1,},](v0)(v1)
\Edge[lw=0.1cm,style={color=cv0v5,},](v0)(v5)
\Edge[lw=0.1cm,style={color=cv1v2,},](v1)(v2)
\Edge[lw=0.1cm,style={color=cv1v4,},](v1)(v4)
\Edge[lw=0.1cm,style={color=cv2v3,},](v2)(v3)
\Edge[lw=0.1cm,style={color=cv5v6,},](v5)(v6)
\Edge[lw=0.1cm,style={color=cv6v7,},](v6)(v7)
\end{tikzpicture}
\caption{Bipartite graph on $8$ vertices}
\label{not return pst graph}
\end{center}
\end{figure}
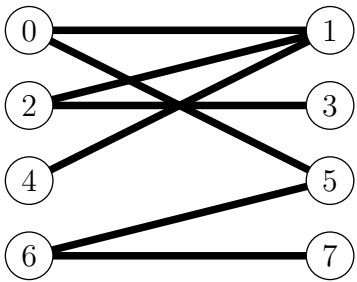

Now consider the bipartite graph $G$ in Figure~\ref{not return pst graph} as an example. We define a bipartite walk on $G$. The two parts of $G$ are $C_0=\{0,2,4,6\}$ and $C_2=\{1,3,5,7\}$. Here is an example of how we partition edges of $G$: edges $(0,1),(0,5)$ are in the same cell in parition $\pi_0$ and edges $(0,1),(2,1),(4,1)$ are in the same cell in partition $\pi_1$.
We have that  
\[
	\widehatP_0=\begin{pmatrix}
 \frac{1}{\sqrt{2}} & 0 & 0 & 0 \\
 \frac{1}{\sqrt{2}}& 0 & 0 & 0 \\
0 & \frac{1}{\sqrt{2}}& 0 & 0 \\
0 & 0 & 1 & 0 \\
0 &  \frac{1}{\sqrt{2}}& 0 & 0 \\
0 & 0 & 0 &  \frac{1}{\sqrt{2}}\\
0 & 0 & 0 &  \frac{1}{\sqrt{2}}
\end{pmatrix},\quad
\widehatP_1=\begin{pmatrix}
\frac{1}{\sqrt{3}} & 0 & 0 & 0 \\
0 & 0 & \frac{1}{\sqrt{2}} & 0 \\
\frac{1}{\sqrt{3}} & 0 & 0 & 0 \\
\frac{1}{\sqrt{3}}& 0 & 0 & 0 \\
0 & 1 & 0 & 0 \\
0 & 0 & \frac{1}{\sqrt{2}}& 0 \\
0 & 0 & 0 & 1
\end{pmatrix} 
\]
and we get the corresponding projections
\[
P=
\begin{pmatrix}
\frac{1}{3}& 0 & \frac{1}{3}& \frac{1}{3}& 0 & 0 & 0 \\[2.5mm]
0 & \frac{1}{2} & 0 & 0 & 0 & \frac{1}{2}& 0 \\[2.5mm]
\frac{1}{3}& 0 & \frac{1}{3} & \frac{1}{3}& 0 & 0 & 0 \\[2.5mm]
\frac{1}{3}& 0 & \frac{1}{3}&\frac{1}{3} & 0 & 0 & 0 \\[2.5mm]
0 & 0 & 0 & 0 & 1 & 0 & 0 \\[2.5mm]
0 & \frac{1}{2}& 0 & 0 & 0 & \frac{1}{2} & 0 \\[2.5mm]
0 & 0 & 0 & 0 & 0 & 0 & 1
\end{pmatrix},\quad
Q=\begin{pmatrix}
\frac{1}{2}& \frac{1}{2} & 0 & 0 & 0 & 0 & 0 \\[2.5mm]
\frac{1}{2}& \frac{1}{2}& 0 & 0 & 0 & 0 & 0 \\[2.5mm]
0 & 0 & \frac{1}{2}& 0 & \frac{1}{2} & 0 & 0 \\[2.5mm]
0 & 0 & 0 & 1 & 0 & 0 & 0 \\[2.5mm]
0 & 0 & \frac{1}{2} & 0 & \frac{1}{2} & 0 & 0 \\[2.5mm]
0 & 0 & 0 & 0 & 0 & \frac{1}{2} & \frac{1}{2} \\[2.5mm]
0 & 0 & 0 & 0 & 0 & \frac{1}{2} & \frac{1}{2}
\end{pmatrix}.
\]
The transition matrix of the bipartite walk on $G$ is 
\[
	U=\begin{pmatrix}
0 & -\frac{1}{3} & 0 &  \frac{2}{3}&  \frac{2}{3} & 0 & 0 \\[2.5mm]
0 & 0 & 0 & 0 & 0 & 0 & 1 \\[2.5mm]
0 & \frac{2}{3}  & 0 & \frac{2}{3}  & \frac{1}{3} & 0 & 0 \\[2.5mm]
0 &  \frac{2}{3}  & 0 & -\frac{1}{3}  &  \frac{2}{3}  & 0 & 0 \\[2.5mm]
0 & 0 & 1 & 0 & 0 & 0 & 0 \\[2.5mm]
1 & 0 & 0 & 0 & 0 & 0 & 0 \\[2.5mm]
0 & 0 & 0 & 0 & 0 & 1 & 0
\end{pmatrix}.
\]

\section{Spectrum of a biregular $G$ determines the spectrum of the bipartite walk on $G$}

Given a graph $G$, we define the bipartite walk over $G$ with the transition matrix $U$. In this section, we show that if $G$ is biregular, the spectrum of the adjacency matrix of $G$ determines the spectrum of $U$.

Following the notations defined in previous section, we define 
\[ 
C=P_1^TP_0
\]
and define $\widehat{C} ={\widehatP_1}^T\widehatP_0.$ First we show that the spectrum of $\widehat{C}\widehat{C}^T$ determines the spectrum of $U$. All the theorems cited in this section are stated in terms of bipartite walks. They are originally stated and proved by Zhan in~\cite{harmonyphd} when she discusses the spectra of matrices that are products of two reflections. Detailed proofs are given in~\cite{harmonyphd}, so we omit them here.

\begin{theorem}[Lemma~$2.3.5$ in~\cite{harmonyphd}]
\label{1-eignsp of U}
 	Let $P,Q$ be projections on $\mathbb{C}^m$.
	The $1$-eigenspace of $U$ is
	 \[ 
		\left(\col(P)\cap \col(Q)\right)\oplus \left(\ker(P)\cap \ker(Q)		\right)
	\] and it has dimension
 \[	
	m-\rk(P)-\rk(Q)+2\dim\left(\col(P)\cap \col(Q)\right).
\]
Moreover, the map $y\mapsto \widehatP_1 y$ is an isomorphism from the $1$-eigenspace of $\widehat{C}\widehat{C}^T$ to $\col(P)\cap \col(Q)$.\qed
\end{theorem}

Recall that $C_0,C_1$ are two color classes of the underlying bipartite graph $G$.
\begin{theorem}[Lemma $2.3.6$ in~\cite{harmonyphd}]
\label{-1 eigen}
	The $(-1)$-eigenspace for $U$  is 
	\[
		\left(\col(P)\cap \ker(Q)\right)\oplus \left(\ker(P)\cap \col(Q)	\right)\] and its dimension is 
		\[
			\abs{C_0}+\abs{C_1}-2\rk(\widehat{C}).
			\]
Moreover, the map $y\mapsto\widehatP_0$ is an isomorphism from $\ker(\widehat{C})$ to $\col(P)\cap \ker(Q)$, and the map $y\mapsto {\widehat{P}^*_1}y$ is an isomorphism from $\ker(\widehat{C}^T)$ to $\ker(P)\cap \col(Q)$.\qed
\end{theorem}

\begin{theorem}[Lemma $2.3.4$ in \cite{harmonyphd}]
\label{e^itheta-eigenspace of U}
Let $\mu\in[0,1]$ be an eigenvalue of $\widehat{C}\widehat{C}^T$. Choose $\theta$ such that $\cos\theta=2\mu-1.$ Let $E_\mu$ be the orthogonal projection onto the $\mu$-eigenspace of $\widehat{C}\widehat{C}^T$. Set 
\[
	W:= \widehatP_1E_\mu\widehatP_1^T. 
\]
Then the $e^{i\theta}$-eigenmatrix of $U$ is
\[
	\frac{1}{\sin^2(\theta)}\left( (\cos\theta+1)W-(e^{i\theta}+1)PW-(e^{-i		\theta}+1)WP+2PWP\right),
\] 
and the $e^{-i\theta}$-eigenmatrix of $U$ is
\[
	\frac{1}{\sin^2(\theta)}\left( (\cos\theta+1)W-(e^{-i\theta}+1)PW-(e^{i	\theta}+1)WP+2PWP\right).\qed
\]
\end{theorem}

We have shown that the spectrum of $\widehat{C}\widehat{C}^T$ determines the spectrum of $U$. Now we are going to show that the spectrum of $A(G)$ determines the spectrum of $\widehat{C}\widehat{C}^T$ when $G$ is a biregular graph.

\begin{theorem} 
\label{spec(A) determines spec(U)}
Let $G$ be a biregular graph with degree $(d_0,d_1)$ and $U$ is the transition matrix of the bipartite walk defined over $G$.Then 
For every complex eigenvalue $e^{i\theta_r}$ of $U$, we have that 
\[
\cos\theta=2\frac{\lambda^2}{d_0d_1}-1,
\] where $\lambda_r$ is an eigenvalue of $A(G)$.
\end{theorem}
\proof 
Since
\[ 
C=P_1^TP_0,
\]
using the definition of $P,Q$, it is not hard to see that
\[
A(G)=
\begin{pmatrix}
\mathbf{0}&C\\
C^T&\mathbf{0}
\end{pmatrix}.
\] 
Recall that $\widehat{C} ={\widehat{P}_1}^T\widehat{P}_0$. If $G$ is biregular with degree $(d_0,d_1)$, we have that \[
\widehat{C}=\frac{1}{\sqrt{d_0d_1}}C.
\]
It follows that
\[
A^2=
\begin{pmatrix}
\mathbf{0}&C\\
C^T&\mathbf{0}
\end{pmatrix} \begin{pmatrix}
\mathbf{0}&C\\
C^T&\mathbf{0}
\end{pmatrix}=
\begin{pmatrix}
CC^T&\mathbf{0}\\
\mathbf{0}&C^TC
\end{pmatrix},
\]
which implies that 
\begin{equation}
\label{spec of A2 and CCT}
\spec(A^2)=\spec(CC^T)\cup\spec(CC^T) .
\end{equation}

Theorem~\ref{e^itheta-eigenspace of U} states that for every eigenvalue $\mu$ of $\widehat{C}\widehat{C}^T$ with $\mu\in(0,1)$ , we have that 
\[
	\cos\theta=2\mu-1
	\] 
for every complex eigenvalue $e^{i\theta}$ of $U$. Since $G$ is biregular with degree $(d_0,d_1)$, we have that 
\[
\frac{1}{d_0d_1}CC^T=\widehat{C}\widehat{C}^T.
\]
For every complex eigenvalue $e^{i\theta}$ of $U$, by relation~\ref{spec of A2 and CCT},  we have that 
\[
\cos\theta=2\frac{\lambda^2}{d_0d_1}-1,
\] where $\lambda$ is an eigenvalue of $A(G)$.\qed

\section{Periodicity of bipartite walks}
Given a graph $G$, let $U$ denote the transition matrix of the bipartite walk defined over $G$. Note that the rows and columns of $U$ are indexed by the edges of $G$. Let $e_a$ denote the standard basis vector in $\mathbb{C}^{E(G)}$ indexed by the edge $a$ in graph $G$. The matrix 
\[
D_a=e_ae_a^T
\] 
is a state associated with edge $a$ of $G$. In the bipartite walk whose transition matrix is $U$, if a walker starts at state $D_a$, after $k$ steps, the walker is at the state 
\[
D_a(k) =U^k D_a {U^*}^k .\]
 We say a state $a$ is \textsl{periodic} if and only if there exist a positive integer $\tau$ such that 
\[
D_a(\tau)=
U(\tau) D_a U(-\tau) =D_a
.\]

The walk over the graph $G$ is periodic if every state of $G$ is periodic, i.e., there exist a positive integer $k$ such that 
\[
U^k=I.\] 
In this case, sometimes we also say the graph is periodic when the context is clear.
Periodicity is an important and interesting property of quantum walks. 

In this section, we study the periodicity of bipartite walks. Note that when we talk about biregular graph in the context of bipartite walk, we assume that all the vertices in the same color class have the same degree. We show that if a biregular graph $G$ has eigenvalues whose squares are algebraic integers with degree at most two, there is a characterization of periodicity of bipartite walks in terms of spectrum of $G$. As we stated before, this result extends the result proved by Kubota in~\cite{ShoPED}. 

The \textsl{eigenvalue support} of a state $D$ is the set 
\[
\big\{(\theta_r,\theta_s): E_rDE_s\neq 0\big\}.
\]
\begin{theorem}
\label{PED characterization}
Let $U=\sum_{r} \theta_rE_r$ be the transition matrix of a bipartite walk. 
State $D_a$ is periodic if and only if $\theta_r,\theta_s\in\mathbb{Q}\pi$ for all $(\theta_r,\theta_s)$ in the eigenvalue support of $D_a$.
\end{theorem}
\proof Since $D_a$ is periodic, there is a positive integer $\tau$ such that 
\[
	D_a(\tau)=\sum_{r,s} e^{i\tau(\theta_r-\theta_s)} E_rD_aE_s=D_a.
	\]
Since 
\[
D_a=\sum_{r,s} E_r D_a E_s,\] for every $\theta_r,\theta_s$,
\[
e^{i\tau(\theta_r-\theta_s)} E_rD_aE_s= E_r D_a E_s.\]
Then we have that 
\[
	E_r D_a \overline{E_r}=  e^{2\tau\theta_r i} E_rD_a\overline{E_r}.
	\] 
Since entries of $E_r D_a \overline{E_r}=(E_re_a)\overline{(E_re_a)}$ are the norm squares of entries of $E_re_a$, they are real.
We must have that 
\[
e^{i\tau(\theta_r-\theta_s)}=1,
\]
which implies that 
\[
\tau(\theta_r-\theta_s)=2m_{r,s}\pi
\] for some integer $m_{r,s}$.
In particular, the pair $(\theta_r,\theta_s)$, where $\theta_s=-\theta_r$, is in the eigenvalue support of $D_a$, so we have that 
\[
2\tau\theta_r=2m_{r,s}\pi,
\] 
for some integer $m_r$. Therefore,
\[
	\theta_r=\frac{m_{r,s}}{\tau}\pi\in\mathbb{Q}\pi.\qed
	\]

From now on, we consider the case when the graph is periodic, i.e. every state of $G$ is periodic. That is, we want that Theorem~\ref{PED characterization} applies to all the eigenvalues of $G$.
%\begin{theorem}[Theorem $1$ in~\cite{trigoometricalgebraicnum}]
%\label{2cos alge int}
%Let $r=\frac{k}{n}$, where $n>2$, be a rational number with $\gcd(k,n)=1$. Let $\phi(n)$ be the number of integers less than $n$ and relatively prime to $n$. Then $2\cos(2\pi r)$ in an algebraic integer with degree $\phi(n)/2$.
%\end{theorem}

The following theorem is essentially the same as Lemma $3.2$ in~\cite{ShoPED}.

\begin{theorem}
\label{2cos alge int}
 Let $U=\sum_r e^{i\theta_r}E_r$ be the transition matrix of the bipartite walk defined over graph $G$. If the bipartite walk is periodic, then $2\cos k\theta_r$ is an algebraic integer for any non-negative integer $k$.
\end{theorem}
\proof
If $G$ is periodic with period $t$, i.e.,
\[
U^t=I,\]
 eigenvalue
$e^{i\theta_r}$ of $U$ is a root of $x^t-1$ and this implies that $e^{i\theta_r}$ is an algebraic integer.
 %Let $\mu_r\in(0,1)$ be an eigenvalue of $\widehat{C}\widehat{C}^T$.
 Then for any non-negative integer $k$, we have that
\[
\left(e^{i\theta_r}\right)^k+\left( e^{-i\theta_r}\right)^k=2\cos k\theta_r\] is an algebraic integer.\qed

Using the theorem above, we can derive a necessary condition for a graph being periodic. Computationally, this provides a easy way for us to determine when a graph is not periodic.

\begin{theorem}
\label{ped trace}
Let $G$ be a periodic bipartite graph. Then $\tr(U^k)\in\mathbb{Z}$ for any integer $k$.
\end{theorem}
\proof Since $U$ is a rational matrix, we know $\tr(U)\in\mathbb{Q}$. On the other hand, 
\[
\tr(U^k)=\sum_r e^{ik\theta_r}=\sum_r 2\cos k\theta_r .\] 
By Theorem~\ref{2cos alge int}, we know
$ 2\cos k\theta_r$ is an algebraic integer for any integer $k$. So $\tr(U^k)$ is an algebraic integer. Hence, when $G$ is periodic, we must have that 
\[
\tr(U^k)\in\mathbb{Z}
\]
for any integer $k$.
\qed

Using the necessary condition for periodicity stated above, we can easily see that the graph shown in Figure~\ref{not return pst graph} is not periodic, since the transition matrix has trace $-\frac{1}{3}$.

%\begin{corollary}
%Let $G$ be a biregular graph with degree $(d_0,d_1)$ and $\lambda$ is an eigenvalue of $A(G)$. The bipartite walk defined over $G$ is periodic, 
% \[
% 2\cos\theta=
%4\frac{\lambda^2}{d_0d_1}-2
%\]
%is an algebraic integer.
%\end{corollary}

We are going to use the following two results from algebraic number theory to give a characterization of periodicity of bipartite walks over biregular graph that squares of their eiganvalues are algebraic integer with degree at most two.

\begin{lemma}[Theorem $3.3$ in~\cite{algtripi}]\label{cos rational pi}
Let $\alpha\in[-1,1]$. Assume that $\frac{1}{\pi}\arccos\alpha=\frac{2k}{n}, k\in\mathbb{Z}, n\in\mathbb{N}, \gcd(k,n)=1$. Then 
\begin{enumerate}[label=(\roman*)]
\item the number $2\alpha=2\cos\frac{2k\pi}{n}$ is an algebraic integer of degree one if and only if $n=1,2,3,4,6$; in such cases, all the values taken by $\alpha$ are 
\begin{align*}
1&=\cos0=-\cos\pi,\quad 0 =\cos\frac{2\pi}{4}=\cos\frac{6\pi}{4},\\
\frac{1}{2}&=\cos\frac{2\pi}{6}=\cos\frac{10\pi}{6}=-\cos\frac{2\pi}{3}=-\cos\frac{4\pi}{3};
\end{align*}
\item the number $2\alpha=2\cos\frac{2k\pi}{n}$ is an algebraic integer of degree two if and only if $n=5,8,10,12$; in such cases, all the values taken by $\alpha$ are 
\begin{align*}
\frac{\sqrt{5}-1}{4}&=\cos\frac{2\pi}{5}=\cos\frac{8\pi}{5}=-\cos\frac{6\pi}{10}=-\cos\frac{14\pi}{10},\\
\frac{\sqrt{5}+1}{4}&=-\cos\frac{4\pi}{5}=-\cos\frac{6\pi}{5}=\cos\frac{2\pi}{10}=\cos\frac{18\pi}{10},\\
\frac{\sqrt{2}}{2}&=\cos\frac{2\pi}{8}=\cos\frac{14\pi}{8}=-\cos\frac{6\pi}{8}=-\cos\frac{10\pi}{8},\\
\frac{\sqrt{3}}{2}&=\cos\frac{2\pi}{12}=\cos\frac{22\pi}{12}=-\cos\frac{10\pi}{12}=-\cos\frac{14\pi}{12}.
\end{align*}\qed
\end{enumerate} 
\end{lemma}

\begin{lemma}[Proposition $2.34$ in~\cite{jarvis2014algebraic}]
\label{form of Q(sqrt(m)) alg int}
Suppose that $m$ is a square-free integer (i.e., not divisible by the square of any prime). Let $\Omega$ denote the set of algebraic integer. Then 
\[
\Omega\cap\mathbb{Q}(\sqrt{m})=
\begin{cases}
\{p+q\sqrt{m}: p,q \in \mathbb{Z}\}\quad\text{if }m\equiv 2,3 \pmod 4,\\[3mm]
\{p+\frac{1+\sqrt{m}}{2}q: p,q \in \mathbb{Z}\}\quad\text{if }m\equiv 1 \pmod 4.
\end{cases}\qed
\]
\end{lemma}

\begin{theorem}
\label{biregular ped}
Let $G$ be a biregular graph with degree $(d_0,d_1)$. Assume that squares of eigenvalues $\lambda_r$ of $A(G)$ are algebraic integers with degree at most two.  Consider bipartite walk defined on $G$. The bipartite walk defined over $G$ is periodic if and only if every eigenvalue $\lambda_r$ of $A(G)$ satisfies that
\begin{enumerate}[label=(\alph*)]
\item if $\lambda_r^2$ is an algebraic integer of degree one, then $d_0d_1\equiv 0\pmod 4$ and
\[
	\lambda_r^2 \in \bigg\{\frac{1}{2}d_0d_1, \frac{3}{4}d_0 d_1,\frac{1}{4} d_0 d_1,0,d_0 d_1\bigg\}
	\];
\item if $\lambda_r^2$ is an algebraic integer of degree two, then 
\[
	\lambda_r^2\in \Bigg\{\left(\frac{1}{2}\pm\frac{\sqrt{2}}{4}\right)d_0d_1,\left(\frac{1}{2}\pm\frac{\sqrt{3}}{4}\right)d_0d_1,
\frac{5\pm\sqrt{5}}{8}d_0d_1 ,\frac{3\pm\sqrt{5}}{8}d_0d_1\Bigg\}.
\]
Moreover, $\lambda_r^2$ comes in algebraic conjugate pairs.
\end{enumerate}
\end{theorem}
\proof
By Theorem~\ref{PED characterization}, the bipartite walk define on $G$ is periodic if and only if $\frac{\theta_r}{\pi}\in\mathbb{Q}$ for all eigenvalues $e^{i\theta_r}$ of $U$. It is easy to check that the condition is sufficient. Here we are going to prove it is necessary.

By Theorem~\ref{2cos alge int}, if $G$ is periodic, then $2\cos\theta_r$ is an algebraic integer.
As shown in Theorem~\ref{spec(A) determines spec(U)}, for every eigenvalue $\lambda_r$ of $A(G)$,
\[
\cos\theta_r=2\frac{\lambda_r^2}{d_0d_1}-1.
\]
By assumption, the number $\lambda_r^2$ is an algebraic integer of degree at most two, which is the same as
\[
\lambda_r^2=a+b\sqrt{m_r}
\] for some square-free integer $m_r$ and $a,b\in \mathbb{Q}$. 
When $\lambda^2$ is either in $\mathbb{Q}$ or $\lambda_r^2=a+b\sqrt{m_r}$ for some square-free integer $m_r$ and non-zero $a,b\in \mathbb{Q}$, by Lemma~\ref{cos rational pi}, we know that $G$ is periodic if and only if \begin{enumerate}[label=(\alph*)]
\item when $b=0$, i.e., $\lambda_r^2 \in \mathbb{Q}$, 
\[
\cos\theta_r \in \big \{ 0,\pm 1,\pm\frac{1}{2}\big\};
\]
\item when $\lambda_r^2=a+b\sqrt{m_r}$ for some square-free integer $m_r$ and non-zero $a,b\in \mathbb{Q}$, 
\[
\cos\theta_r \in \big\{\pm \frac{\sqrt{2}}{2},\pm\frac{\sqrt{3}}{2},\pm\frac{\sqrt{5}+1}{4},\pm \frac{\sqrt{5}-1}{4}\big\}
.\]
\end{enumerate}

First consider the case when $\lambda_r^2\in \mathbb{Q}$.
When 
\[
\cos\theta=2\frac{\lambda^2}{d_0d_1}-1=0,
\]
 we have that
 \[
 \lambda^2=\frac{1}{2}d_0d_1.\]
  When 
\[
\cos\theta=2\frac{\lambda^2}{d_0d_1}-1=1,
\]
 we have that
 \[
\lambda^2={d_0d_1},\] which is guaranteed since the largest eigenvalue of $A(G)$ is $\sqrt{d_0d_1}$. Similarly, when $\cos\theta=2\frac{\lambda^2}{d_0d_1}-1=-1$, we have that $\lambda_r=0$.

Consider the case when 
\[\cos\theta=
2\frac{\lambda_r^2}{d_0d_1}-1=\frac{1}{2},
\] and
 we have that
 \[
 \frac{\lambda_r^2}{d_0d_1}=\frac{3}{4}.\]Similarly, when $\cos\theta=-\frac{1}{2}$, we have that $\lambda_r^2=\frac{1}{4}d_0d_1$.

Now consider the case when $\lambda_r^2=a+b\sqrt{m_r}$ for some square-free integer $m_r$ and non-zero $a,b\in \mathbb{Q}$.
Assume
\[
\cos\theta=
2\left(\frac{\lambda_r^2}{d_0d_1}-1\right)=\pm\frac{\sqrt{2}}{2}.
\] Using Lemma~\ref{form of Q(sqrt(m)) alg int}, we assume that
\[
\lambda_r^2=p+q\sqrt{2},\] where $p,q$ are both non-zero integers
and 
\[
\cos\theta=2\frac{(p+q\sqrt{2})}{d_0d_1}-1=\pm\frac{\sqrt{2}}{2}
\]
Then we have that \[
\frac{2p}{d_0d_1}-1=0,\quad \frac{2q}{d_0d_1}=\pm\frac{1}{2}.
\]
Combining both equations above, we have that 
\[
p=\pm 2q= d_0d_1
\]
and \[
\lambda^2=(1\pm\frac{\sqrt{2}}{2})p=\left(\frac{1}{2}\pm\frac{\sqrt{2}}{4}\right)d_0d_1.\]

When
\[
2\frac{\lambda_r^2}{d_0d_1}-1=\pm\frac{\sqrt{3}}{2},\] 
by Lemma~\ref{form of Q(sqrt(m)) alg int}, we assume that
\[
\lambda_r^2 = p+q\sqrt{3}.
\] Then using the similar argument as previous case, we have that
\[
\lambda^2=(1\pm \frac{\sqrt{3}}{2})p=\left(\frac{1}{2}\pm\frac{\sqrt{3}}{4}\right)d_0d_1.\]

Now, consider the case when
\[
\cos\theta_r=2\frac{\lambda_r^2}{d_0d_1}-1=\pm\frac{\sqrt{5}+1}{4}.
\]
This implies that
\[
\lambda_r^2=p+\frac{1+\sqrt{5}}{2}q, 
\]
where $p,q$ are both non-zero integers.
So
\[
\frac{2}{d_0d_1} \left(p+\frac{1+\sqrt{5}}{2}q\right)-1=\pm\frac{1+\sqrt{5}}{4}
\]
This implies that
\[
\frac{2}{d_0d_1}(p+\frac{q}{2})-1=\frac{1}{4},\quad
\frac{q}{d_0d_1}=\frac{1}{4}.
\]
or
\[
\frac{2}{d_0d_1}(p+\frac{q}{2})-1=-\frac{1}{4},\quad
\frac{q}{d_0d_1}=-\frac{1}{4}.
\]
Combining two equation above, we have
\[
p=\pm 2q=\frac{1}{2}d_0d_1,
\]
and consequently,\[
\lambda_r^2=\frac{5+\sqrt{5}}{8}d_0d_1 \text{ or }\frac{3-\sqrt{5}}{8}d_0d_1
\] when $p=2q$ and $p=-2q$ respectively. Similarly, when $\cos\theta_r=\pm\frac{\sqrt{5}-1}{4}$, we have 
\[
p=-3q=\frac{3}{4}d_0d_1,\text{ or } p=q=\frac{1}{4}d_0d_1,
\]
and consequently,
\[
\lambda^2= \frac{5-\sqrt{5}}{8}d_0d_1, \text{ or } \lambda^2=\frac{3+\sqrt{5}}{8}d_0d_1.\]

We can view $\lambda^2$ as eigenvalues of $A(G)^2$. So $\lambda^2$ comes in algebraic conjugate pairs.
\qed

\section{Grover's walk is a special case of bipartite walk}
\label{GW subset BW}

Grover's walk is a well-studied discrete quantum walk model. We are going to show that Grover's walk is a special case of bipartite walk model. That is, given a graph $G$, the transition matrix of the Grover's walk on $G$ is the same as the transition matrix of the bipartite walk on the subdivision graph of $G$.

First, we define the transition matrix of Grover's walks. Given a graph $G$, we can give directions to the edges of $G$ such that the arc set of $G$ is 
\[
\mathcal{A}=\big\{ (a,b), (b,a) \mid \{a,b\}\in E(G)\big\}.
\] Let $\alpha=(x,y)$ be an arc in $\mathcal{A}$, we say $x$ is the head of $\alpha$, denoted by $o(\alpha)$ and $y$ is the tail of $\alpha$, denoted by $t(\alpha)$. When $\alpha=(x,y)$, we say $\alpha^{-1}=(y,x)$.
Define a matrix $D\in \mathbb{C}^{V\times \mathcal{A}}$ such that 
\[
D_{x,\alpha}=\frac{1}{\sqrt{\deg(x)}}\delta_{x,t(\alpha)}.
\] Then $D^*D\in\mathbb{C}^{\mathcal{A}\times \mathcal{A}}$ are
\[
\left(D^*D\right)_{\alpha,\beta}
=\begin{cases}
\frac{1}{\deg\left(t(\alpha)\right)} \quad \text{if }t(\alpha)=t(\beta)\\[3mm]
0 \quad\text{otherwise.}
\end{cases}
\]
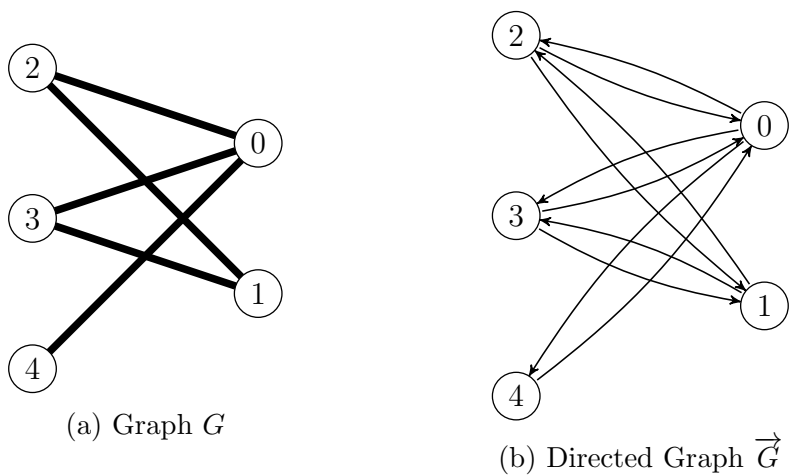
\begin{figure}[H]
\centering
\begin{subfigure}{0.4\linewidth}
\centering
\begin{tikzpicture}[scale=1]
\definecolor{cv0}{rgb}{0,0,0}
\definecolor{cfv0}{rgb}{1,1,1}
\definecolor{clv0}{rgb}{0,0,0}
\definecolor{cv1}{rgb}{0,0,0}
\definecolor{cfv1}{rgb}{1,1,1}
\definecolor{clv1}{rgb}{0,0,0}
\definecolor{cv2}{rgb}{0,0,0}
\definecolor{cfv2}{rgb}{1,1,1}
\definecolor{clv2}{rgb}{0,0,0}
\definecolor{cv3}{rgb}{0,0,0}
\definecolor{cfv3}{rgb}{1,1,1}
\definecolor{clv3}{rgb}{0,0,0}
\definecolor{cv4}{rgb}{0,0,0}
\definecolor{cfv4}{rgb}{1,1,1}
\definecolor{clv4}{rgb}{0,0,0}
\definecolor{cv0v2}{rgb}{0,0,0}
\definecolor{cv0v3}{rgb}{0,0,0}
\definecolor{cv0v4}{rgb}{0,0,0}
\definecolor{cv1v2}{rgb}{0,0,0}
\definecolor{cv1v3}{rgb}{0,0,0}
\Vertex[style={minimum size=1cm,draw=cv0,fill=cfv0,text=clv0,shape=circle},LabelOut=false,L=\hbox{$0$},x=3cm,y=4cm]{v0}
\Vertex[style={minimum size=1cm,draw=cv1,fill=cfv1,text=clv1,shape=circle},LabelOut=false,L=\hbox{$1$},x=3cm,y=2.0cm]{v1}
\Vertex[style={minimum size=1cm,draw=cv2,fill=cfv2,text=clv2,shape=circle},LabelOut=false,L=\hbox{$2$},x=0cm,y=5cm]{v2}
\Vertex[style={minimum size=1cm,draw=cv3,fill=cfv3,text=clv3,shape=circle},LabelOut=false,L=\hbox{$3$},x=0cm,y=3cm]{v3}
\Vertex[style={minimum size=1cm,draw=cv4,fill=cfv4,text=clv4,shape=circle},LabelOut=false,L=\hbox{$4$},x=0cm,y=1cm]{v4}
\Edge[lw=0.1cm,style={color=cv0v2,},](v0)(v2)
\Edge[lw=0.1cm,style={color=cv0v3,},](v0)(v3)
\Edge[lw=0.1cm,style={color=cv0v4,},](v0)(v4)
\Edge[lw=0.1cm,style={color=cv1v2,},](v1)(v2)
\Edge[lw=0.1cm,style={color=cv1v3,},](v1)(v3)
\end{tikzpicture}
\caption{ Graph $G$}
\label{G}
\end{subfigure}
~
\begin{subfigure}{0.5\linewidth}
\centering
\begin{tikzpicture}[scale=0.6]
\definecolor{cv0}{rgb}{0,0,0}
\definecolor{cfv0}{rgb}{1,1,1}
\definecolor{clv0}{rgb}{0,0,0}
\definecolor{cv1}{rgb}{0,0,0}
\definecolor{cfv1}{rgb}{1,1,1}
\definecolor{clv1}{rgb}{0,0,0}
\definecolor{cv2}{rgb}{0,0,0}
\definecolor{cfv2}{rgb}{1,1,1}
\definecolor{clv2}{rgb}{0,0,0}
\definecolor{cv3}{rgb}{0,0,0}
\definecolor{cfv3}{rgb}{1,1,1}
\definecolor{clv3}{rgb}{0,0,0}
\definecolor{cv4}{rgb}{0,0,0}
\definecolor{cfv4}{rgb}{1,1,1}
\definecolor{clv4}{rgb}{0,0,0}
\definecolor{cv0v2}{rgb}{0,0,0}
\definecolor{cv0v3}{rgb}{0,0,0}
\definecolor{cv0v4}{rgb}{0,0,0}
\definecolor{cv1v2}{rgb}{0,0,0}
\definecolor{cv1v3}{rgb}{0,0,0}
\definecolor{cv2v0}{rgb}{0,0,0}
\definecolor{cv2v1}{rgb}{0,0,0}
\definecolor{cv3v0}{rgb}{0,0,0}
\definecolor{cv3v1}{rgb}{0,0,0}
\definecolor{cv4v0}{rgb}{0,0,0}
\Vertex[style={minimum size=1cm,draw=cv0,fill=cfv0,text=clv0,shape=circle},LabelOut=false,L=\hbox{$0$},x=5.5cm,y=5cm]{v0}
\Vertex[style={minimum size=1cm,draw=cv1,fill=cfv1,text=clv1,shape=circle},LabelOut=false,L=\hbox{$1$},x=5.5cm,y=1cm]{v1}
\Vertex[style={minimum size=1cm,draw=cv2,fill=cfv2,text=clv2,shape=circle},LabelOut=false,L=\hbox{$2$},x=0cm,y=7cm]{v2}
\Vertex[style={minimum size=1cm,draw=cv3,fill=cfv3,text=clv3,shape=circle},LabelOut=false,L=\hbox{$3$},x=0cm,y=3cm]{v3}
\Vertex[style={minimum size=1cm,draw=cv4,fill=cfv4,text=clv4,shape=circle},LabelOut=false,L=\hbox{$4$},x=0cm,y=-1cm]{v4}
\Edge[lw=0.1cm,style={post, bend right=8,color=cv0v2,},](v0)(v2)
\Edge[lw=0.1cm,style={post, bend right=10,color=cv0v3,},](v0)(v3)
\Edge[lw=0.1cm,style={post, bend right=10,color=cv0v4,},](v0)(v4)
\Edge[lw=0.1cm,style={post, bend right=8,color=cv1v2,},](v1)(v2)
\Edge[lw=0.1cm,style={post, bend right=10,color=cv1v3,},](v1)(v3)
\Edge[lw=0.1cm,style={post, bend right=8,color=cv2v0,},](v2)(v0)
\Edge[lw=0.1cm,style={post, bend right=8,color=cv2v1,},](v2)(v1)
\Edge[lw=0.1cm,style={post, bend right=10,color=cv3v0,},](v3)(v0)
\Edge[lw=0.1cm,style={post, bend right=10,color=cv3v1,},](v3)(v1)
\Edge[lw=0.1cm,style={post, bend right=10,color=cv4v0,},](v4)(v0)
\end{tikzpicture}
\caption{Directed Graph $\overrightarrow{G}$}
\label{directed G}
\end{subfigure}
\end{figure}
Let $R\in \mathbb{C}^{\mathcal{A}\times \mathcal{A}}$ denote the arc-reversal matrix, i.e.,
\[
R_{\alpha,\beta}=\delta_{\alpha,\beta^{-1}}.
\]
The transition matrix of the Grover's walk defined on $G$ is 
\[
R(2D^*D-I).\]

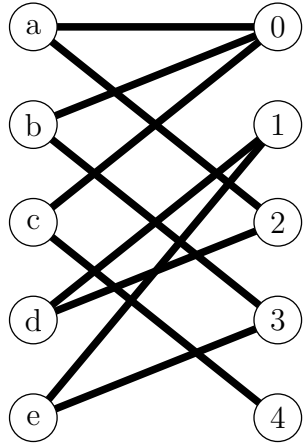
\begin{figure}[H]
\centering
 \begin{tikzpicture}[scale=0.65]
\definecolor{cv0}{rgb}{0,0,0}
\definecolor{cfv0}{rgb}{1,1,1}
\definecolor{clv0}{rgb}{0,0,0}
\definecolor{cv1}{rgb}{0,0,0}
\definecolor{cfv1}{rgb}{1,1,1}
\definecolor{clv1}{rgb}{0,0,0}
\definecolor{cv2}{rgb}{0,0,0}
\definecolor{cfv2}{rgb}{1,1,1}
\definecolor{clv2}{rgb}{0,0,0}
\definecolor{cv3}{rgb}{0,0,0}
\definecolor{cfv3}{rgb}{1,1,1}
\definecolor{clv3}{rgb}{0,0,0}
\definecolor{cv4}{rgb}{0,0,0}
\definecolor{cfv4}{rgb}{1,1,1}
\definecolor{clv4}{rgb}{0,0,0}
\definecolor{cv5}{rgb}{0,0,0}
\definecolor{cfv5}{rgb}{1,1,1}
\definecolor{clv5}{rgb}{0,0,0}
\definecolor{cv6}{rgb}{0,0,0}
\definecolor{cfv6}{rgb}{1,1,1}
\definecolor{clv6}{rgb}{0,0,0}
\definecolor{cv7}{rgb}{0,0,0}
\definecolor{cfv7}{rgb}{1,1,1}
\definecolor{clv7}{rgb}{0,0,0}
\definecolor{cv8}{rgb}{0,0,0}
\definecolor{cfv8}{rgb}{1,1,1}
\definecolor{clv8}{rgb}{0,0,0}
\definecolor{cv9}{rgb}{0,0,0}
\definecolor{cfv9}{rgb}{1,1,1}
\definecolor{clv9}{rgb}{0,0,0}
\definecolor{cv0v5}{rgb}{0,0,0}
\definecolor{cv0v6}{rgb}{0,0,0}
\definecolor{cv0v7}{rgb}{0,0,0}
\definecolor{cv1v8}{rgb}{0,0,0}
\definecolor{cv1v9}{rgb}{0,0,0}
\definecolor{cv2v5}{rgb}{0,0,0}
\definecolor{cv2v8}{rgb}{0,0,0}
\definecolor{cv3v6}{rgb}{0,0,0}
\definecolor{cv3v9}{rgb}{0,0,0}
\definecolor{cv4v7}{rgb}{0,0,0}
\Vertex[style={minimum size=1cm,draw=cv0,fill=cfv0,text=clv0,shape=circle},LabelOut=false,L=\hbox{$0$},x=4cm,y=8cm]{v0}
\Vertex[style={minimum size=1cm,draw=cv1,fill=cfv1,text=clv1,shape=circle},LabelOut=false,L=\hbox{$1$},x=4cm,y=6.0cm]{v1}
\Vertex[style={minimum size=1cm,draw=cv2,fill=cfv2,text=clv2,shape=circle},LabelOut=false,L=\hbox{$2$},x=4cm,y=4cm]{v2}
\Vertex[style={minimum size=1cm,draw=cv3,fill=cfv3,text=clv3,shape=circle},LabelOut=false,L=\hbox{$3$},x=4cm,y=2cm]{v3}
\Vertex[style={minimum size=1cm,draw=cv4,fill=cfv4,text=clv4,shape=circle},LabelOut=false,L=\hbox{$4$},x=4cm,y=0cm]{v4}
\Vertex[style={minimum size=1cm,draw=cv5,fill=cfv5,text=clv5,shape=circle},LabelOut=false,L=\hbox{a},x=-1cm,y=8cm]{v5}
\Vertex[style={minimum size=1cm,draw=cv6,fill=cfv6,text=clv6,shape=circle},LabelOut=false,L=\hbox{b},x=-1cm,y=6cm]{v6}
\Vertex[style={minimum size=1cm,draw=cv7,fill=cfv7,text=clv7,shape=circle},LabelOut=false,L=\hbox{c},x=-1cm,y=4cm]{v7}
\Vertex[style={minimum size=1cm,draw=cv8,fill=cfv8,text=clv8,shape=circle},LabelOut=false,L=\hbox{d},x=-1cm,y=2cm]{v8}
\Vertex[style={minimum size=1cm,draw=cv9,fill=cfv9,text=clv9,shape=circle},LabelOut=false,L=\hbox{e},x=-1cm,y=0cm]{v9}
\Edge[lw=0.1cm,style={color=cv0v5,},](v0)(v5)
\Edge[lw=0.1cm,style={color=cv0v6,},](v0)(v6)
\Edge[lw=0.1cm,style={color=cv0v7,},](v0)(v7)
\Edge[lw=0.1cm,style={color=cv1v8,},](v1)(v8)
\Edge[lw=0.1cm,style={color=cv1v9,},](v1)(v9)
\Edge[lw=0.1cm,style={color=cv2v5,},](v2)(v5)
\Edge[lw=0.1cm,style={color=cv2v8,},](v2)(v8)
\Edge[lw=0.1cm,style={color=cv3v6,},](v3)(v6)
\Edge[lw=0.1cm,style={color=cv3v9,},](v3)(v9)
\Edge[lw=0.1cm,style={color=cv4v7,},](v4)(v7)
\end{tikzpicture}
\caption{The subdivision graph of $G$}
 \label{sub G}
\end{figure}

Given a graph $G$, we define a new graph by subdividing every edge of $G$ exactly once and we call the resulting graph \textsl{the subdivision graph} of $G$, denoted by $S(G)$. For example, the graph in Figure~\ref{sub G} is the subdivision graph of the graph in Figure~\ref{G}.
 Now we are going to show the transition matrix of the Grover's walk defined on $G$ is exactly the same as the transition matrix of the bipartite walk defined on the subdivision graph of $G$.

Given a graph $G$, its subdivision graph $S(G)$ is a bipartite graph with parts $C_0\left(S(G)\right)=V(S(G))\backslash V(G),C_1\left(S(G)\right)=V(G)$. Following the notations and the construction in Section~\ref{BW}, the transition matrix of the bipartite walk on $S(G)$ is \[
U_{S(G)}=(2P_{S(G)}-I)(2Q_{S(G)}-I).\]
Rows and columns of $U_{S(G)}$ are indexed by the edges of $S(G)$.  

Note that for each vertex $i\in C_1\left(S(G)\right)$, we have
$\deg_{S(G)}(i)=\deg_{G}(i)$ and $C_1\left(S(G)\right)=V(G)$.
 Then we index edges of $S(G)$ such that \[
Q_{S(G)}=D^*D.
\] Also, note that for each vertex $a'\in C_0\left(S(G)\right)$, we have
$ \deg_{S(G)}(a')=2$. Two edges $(v_i,a'), (a',v_j)$ of $S(G)$ share a vertex in $C_0\left(S(G)\right)$ if and only if $(v_i,v_j)\in E(G)$. Based on how we index the edge set of $S(G)$, we have
\[
2P_{S(G)}-I=R.\] Thus, we have proved the following theorem.
\begin{theorem}
Given a graph $G$, let $U_{\scriptscriptstyle BW}\left( S(G)\right)$ be the transition matrix of the bipartite walk defined on $S(G)$ and  let $ U_{\scriptscriptstyle GW}\left( G\right)$ be the transition matrix of
Grover's walk defined on $G$.
 Then
\[
U_{\scriptscriptstyle BW}\left( S(G)\right)= U_{\scriptscriptstyle GW}\left( G\right).\qed
\]

\end{theorem}

\section{One step of the biparite walk on $G$ is two steps of Grover's walk on $G$}
\label{1BW=2GW}

Given a bipartite graph $G$, bipartite walks and Grover's walk are two models that one can apply to $G$ to define a discrete quantum walk over $G$. Here we are going to see that two models are closely related, i.e. every $2k$-th power of the transition matrix of Grover's walk on $G$ is a direct sum of $k$-th power of the transition matrix of the bipartite walk defined on the same graph. As shown in the last section, the transition matrix of Grover's walk defined on $G$ equals to the transition matrix of the bipartite walk defined on $S(G)$. So in this section, we are going to look into the connections between the bipartite walk defined on $G$ and the one defined on $S(G)$.

Note that every edges of $S(G)$ has one of its end in $V(G)$. We index rows and columns of  $U_{\scriptscriptstyle BW}\left( S(G)\right)$ such that the first half of the the rows and columns are indexed by the edges with one end in the color class $C_1$ of $G$ and the others are indexed by edges with one end in the color class $C_0$ of $G$. For example, consider the bipartite walk defined on the subdivision graph $S(G)$ shown in Figure~\ref{sub G}. We can index the rows and columns of $Q_{S(G)}$ such that $Q_{S(G)}$ has the form as below.
\[
 \scalemath{0.85}{
\begin{blockarray}{ccccccccccc}
    &(a,0) & (b,0) & (c,0)& (d,1) & (e,1) & (a,2) &(b,3)&(c,4)&(d,2)& (e,3)\\
   \begin{block}{c(ccccc|ccccc)}
(a,0)&   \frac{1}{3} & \frac{1}{3} & \frac{1}{3} & 0 & 0 & 0 & 0 & 0 & 0 & 0 \\[1.5mm]
(b,0)& \frac{1}{3} & \frac{1}{3} & \frac{1}{3} & 0 & 0 & 0 & 0 & 0 & 0 & 0  \\[1.5mm]
(c,0)& \frac{1}{3} & \frac{1}{3} & \frac{1}{3} & 0 & 0 & 0 & 0 & 0 & 0 & 0\\[1.5mm]
 (d,1)& 0 & 0 & 0 & \frac{1}{2} & \frac{1}{2} & 0 & 0 & 0 & 0 & 0  \\[1.5mm]
(e,1)& 0 & 0 & 0 & \frac{1}{2} & \frac{1}{2} & 0 & 0 & 0 & 0 & 0  \\[1.5mm]
\cline{2-11}
(a,2)& 0 & 0 & 0 & 0 & 0 & \frac{1}{2} & 0 & 0 & \frac{1}{2} & 0  \\[1.5mm]
(b,3)& 0 & 0 & 0 & 0 & 0 & 0 & \frac{1}{2} & 0 & 0 & \frac{1}{2} \\[1.5mm]
(c,4)& 0 & 0 & 0 & 0 & 0 & 0 & 0 & 1 & 0 & 0 \\[1.5mm]
(d,2)& 0 & 0 & 0 & 0 & 0 & \frac{1}{2} & 0 & 0 & \frac{1}{2} & 0\\[1.5mm]
(e,3)& 0 & 0 & 0 & 0 & 0 & 0 & \frac{1}{2} & 0 & 0 & \frac{1}{2} \\[1.5mm]
\end{block}
\end{blockarray}}
\]
Let $P_G$ and $Q_G$ denote the projections on two color classes of $G$ as we defined in Section~\ref{BW}. 
Then we can have that 
\[
Q_{S(G)} =\begin{pmatrix}
Q_G& 0\\
0&P_G
\end{pmatrix}.
\]
The transition matrix of bipartite walk defined on $G$ is $U_G$.
Consequently, we have that 
\[ U_{\scriptscriptstyle BW}\left( S(G)\right)= U_{\scriptscriptstyle GW}\left( G\right)=
R\begin{pmatrix}
2Q_G-I& 0\\
0&2P_G-I
\end{pmatrix}=
\begin{pmatrix}
0& 2P_G-I\\
2Q_G-I&0
\end{pmatrix}.
\]

Following directly from the equation above, we have the following theorem.
\begin{theorem}\label{2 step GW is 1 step BW}
Let $G$ be a bipartite graph. Let $U_{\scriptscriptstyle BW}$ and $U_{\scriptscriptstyle GW}$ be the transition matrices of the bipartite walk and Grover's walk defined on $G$ respectively. Then for any non-negative integer $k$, we have that 
\[
{U_{\scriptscriptstyle GW}}^{2k}
=
\begin{pmatrix}
U_{\scriptscriptstyle BW}^k& 0\\[5mm]
0&{U_{\scriptscriptstyle BW}^T}^k
\end{pmatrix}.\qed\]
\end{theorem}

The following corollary shows that the period of a periodic Grover's walk has to be even and the period of a periodic bipartite walk on $G$ is half of the period of the periodic Grover's walk defined on the same graph.
\begin{corollary}
Let $G$ be a bipartite graph. The bipartite walk defined on $G$ is periodic with $\tau$ if and only if Grover's walk defined on $G$ is periodic with period $2\tau$ for some integer $\tau$.
\end{corollary}
\proof 
Let $U_{\scriptscriptstyle BW}$ and $U_{\scriptscriptstyle GW}$ be the transition matrices of the bipartite walk and Grover's walk defined on $G$ respectively. By Theorem~\ref{2 step GW is 1 step BW}, we can see that for any non-negative integer $m$,
\[
U_{\scriptscriptstyle GW}^{2m}=\begin{pmatrix}
U_{\scriptscriptstyle BW}^m& 0\\[5mm]
0&\left(U_{\scriptscriptstyle BW}^T\right)^m
\end{pmatrix}\] and for odd positive integer $i$, we have 
\[
U_{\scriptscriptstyle GW}^{i}=\begin{pmatrix}
0& (2P-I)U_{\scriptscriptstyle BW}^{i-1}\\[5mm]
(2Q-I)\left(U_{\scriptscriptstyle BW}^T\right)^{i-1}&0
\end{pmatrix}.
\]
It follows immediately that Grover's walk on $G$ cannot have odd period.

If the bipartite walk defined on $G$ is periodic with period $\tau$, then it is easy to see that Grover's walk is periodic. Assume Grover's walk on $G$ has period $2k$ for some integer $k$. If $2k<2\tau$, then \[
U_{\scriptscriptstyle BW}^k=I,\]
which contradict that $U_{\scriptscriptstyle BW}$ has period $\tau$. Thus, Grover's walk on $G$ is periodic with period $2\tau$.

Now assume that the Grover's walk is periodic with period $2k$.  By Theorem~\ref{2 step GW is 1 step BW}, the bipartite walk is periodic with period $k$. \qed

\section{Spectrum of $G$ determines the spectrum of $S(G)$}

Let $G$ be a regular graph. Now we are going to give an explicit formula for the eigenvalues of $A(S(G))$ in terms of eigenvalues of $A(G)$. Once we have the formula, by the discussion in Section~\ref{GW subset BW}, we can use Theorem~\ref{biregular ped} to give a characterization of periodicity of Grover's walk on regular graphs whose eigenvalues $\lambda$ satisfies that $\lambda^2$ is a algebraic integer with degree at most two. This characterization extends the result in~\cite{ShoPED}.

Let $B$ be the vertex-edge incidence matrix of $G$ and $\Delta$ is the degree matrix of $G$ . Let $L(G)$ denote the line graph of $G$.
We have that 
\[
BB^T=\left(\Delta(G)+A(G)\right),\quad B^TB=A(L(G))+2I.
\]

The adjacency matrix of the subdivision graph $S(G)$ of $G$ is
\[
\begin{pmatrix}
0& B\\
B^T&0
\end{pmatrix}.
\]
One can easily check that
$\begin{pmatrix}
x&
y
\end{pmatrix}^T$ is an eigenvector of $A(S(G))$ if and only if \[
BB^T x=\lambda x,\quad B^TB y=\lambda y ,
\]
for some $\lambda$. Let $x$ be an eigenvector of $B^TB$ such that $B^TB x=\lambda x$, then \[
\begin{pmatrix}
0& B\\
B^T&0
\end{pmatrix}
\begin{pmatrix}
 Bx\\
\sqrt{\lambda} x
\end{pmatrix}
=
\begin{pmatrix}
\sqrt{\lambda} B x\\
\lambda x
\end{pmatrix}
=
\sqrt{\lambda}
\begin{pmatrix}
 Bx\\
\sqrt{\lambda} x
\end{pmatrix},
\]
which implies that
$\begin{pmatrix}
Bx\\
\sqrt{\lambda} x
\end{pmatrix}$ is an eigenvector of $
\begin{pmatrix}
0& B\\
B^T&0
\end{pmatrix}.$

\begin{lemma}\label{line grp and subdiv grp spectrum}
The value $\lambda-2$ is an eigenvalue of the line graph of $G$ if and only if  $\pm\sqrt{\lambda}$ are eigenvalues of $S(G)$.
\end{lemma}

\proof Eigenvectors of $S(G)$ determines the eigenvectors of $L(G)$. The explicit formula for the correspondence between eigenvalues of $L(G)$ and $S(G)$ derives from the discussion above.\qed

The following lemma is a standard result from algebraic graph theory and this helps us to derive the relation between eigenvalues of $A(G)$ and $A(S(G))$ when $G$ is a regular graph.

\begin{lemma}[Lemma $8.2.5$ in~\cite{Godsil2001}]
\label{evals of LG of regular graph}
Let $G$ be a regular graph of valency $d$ with $n$ vertices and $e$ edges and let $LG$ be the line graph of $G$. Then 
\[
\phi(LG,x)=(x+2)^{e-n}\phi(G,x-d+2).\qed
\]
\end{lemma} 

Since $-2$ is always an eigenvalue of $LG$, zero is always an eigenvalue of $S(G)$.

\begin{corollary}\label{grp and subdiv grp spectrum}
Let $G$ be a $d$-regular graph. Then $\lambda$ is an eigenvalue of $G$ with $\lambda\neq -d$ if and only if $\sqrt{\lambda+d}$ is a non-zero eigenvalue of $S(G)$.
\end{corollary}

\proof Lemma~\ref{evals of LG of regular graph} states that if $\lambda$ is an eigenvalue of $A(G)$, then $\lambda+d-2$ is an eigenvalue of $L(G)$. It follows from Lemma~\ref{line grp and subdiv grp spectrum} that eigenvalue $\lambda$ of $G$ eigenvalue gives the eigenvalue $\sqrt{\lambda+d}$ of $S(G)$.\qed

\section{Periodicity of Grover's walk}
In~\cite{ShoPED}, Kubota studies the periodicity of Grover's walk on regular bipartite graphs with at most five eigenvalues. The eigenvalues of graphs of interest are algebraic integer with degree at most two. Kubota proves the following result.
\begin{theorem}[Theorem $3.3$~\cite{ShoPED}]
\label{SHOresult}
Let $\Gamma$ be a bipartite $k$-regular graph with the $A$-spectrum $\{[\pm k]^{1},[\pm\theta]^{a},[0]^b\}$, where $a\geq 1$ and $b\geq 0$. Then $\Gamma$ is periodic with respect to Grover's walk if and only if  $k$ is even and $\theta\in \{\frac{\sqrt{2}}{2}k,\frac{\sqrt{3}}{2}k\}$.
\end{theorem}

 As we show in Section~\ref{GW subset BW}, given a graph $G$, the transition matrix of the Grover's walk defined over $G$ is the same as the transition matrix of the bipartite walk define over the division graph of $G$. So we can use Theorem~\ref{biregular ped} to extend Theorem~\ref{SHOresult} proved by Kubota. That is, we do not require the underlying graph to be bipartite and we also remove the constraints on the number of distinct eigenvalues the graph can have. We give a characterization of periodicity of Grover's walk on regular graphs whose eigenvalues are algebraic integer with degree at most two. Note that the result in this section are stated in terms of bipartite walk on a subdivision graph of a regular graph.

\begin{corollary} 
\label{subdiv ped}
Let $G$ be a $d$-regular graph, all of whose eigenvalues are algebraic integers of degree at most two in the form of \[
\lambda_r=a+b\sqrt{m_r}\]
for some $a,b\in\mathbb{Q}$ and square-free integer $m_r$. Let $S(G)$ denote the subdivision graph of $G$. The bipartite walk defined over $S(G)$ is periodic if and only if for every eigenvalue $\lambda_r$ of $G$,
\begin{enumerate}[label=(\alph*)]
\item if $b=0$,
$\lambda_r \in \big\{0, \pm d, \pm \frac{1}{2} d\big\}$;
\item if $b\neq 0$, $\lambda_r\in \big\{\pm\frac{\sqrt{2}}{2}d,\pm\frac{\sqrt{3}}{2}d , \frac{1\pm\sqrt{5}}{4}d,\frac{-1\pm \sqrt{5}}{4}d\big\}$.
\end{enumerate}
Note that eigenvalues of $G$ come in algebraic conjugate pairs.
\end{corollary}
\proof By assumption, $\lambda_G$, an eigenvalue of $G$,  is an algebraic integer of degree at most two,  so the eigenvalue $\lambda_{S(G)}$ of $S(G)$ satisfies that  
\[
\lambda_{S(G)}^2=\lambda_G+d.
\] Thus, we know that $\lambda_{S(G)}^2$ is an algebraic integer with degree at most two. We can write 
\[
\lambda_{S(G)}^2=a+b\sqrt{m_r}
\] for some square-free integer $m_r$ and $a,b\in \mathbb{Q}$. 

Graph $S(G)$ is biregular with degree $(2,d)$. By Theorem~\ref{PED characterization}, bipartite walk defined over $S(G)$ is periodic if and only if for every eigenvalue $\lambda_{S(G)}$ of $S(G)$, it satisfies that
\begin{enumerate}[label=(\alph*)]
\item when $b=0$,
$\lambda_{S(G)}^2 \in \big\{d, \frac{3}{2}d,\frac{1}{2} d,0,2d\big\}$;
\item when $b\neq 0$, $\lambda_{S(G)}^2\in  \Big\{\left(1\pm\frac{\sqrt{2}}{2}\right)d,\left(1\pm\frac{\sqrt{3}}{2}\right)d,\frac{\sqrt{5}\pm 5}{4}d,\frac{3\pm\sqrt{5}}{4}d \Big\}$.
\end{enumerate}
Using Corollary~\ref{grp and subdiv grp spectrum}, we see that 
when \[
\lambda_{S(G)}^2=\frac{\sqrt{5}+5}{4}d \text{ or } \frac{3-\sqrt{5}}{4}d,
\]
we get that 
\[
 	\lambda_G=\frac{1+\sqrt{5}}{4}d\text{ or }-\frac{1+\sqrt{5}}{4}d
 	\]
  respectively. Similarly, when 
  \[
  \lambda_{S(G)}^2=\frac{5-\sqrt{5}}{4}d\text{ or }\frac{3+\sqrt{5}}{4}d,
  \]
  we get 
  \[
  \lambda_G=-\frac{\sqrt{5}-1}{4}d\text{ or }\frac{\sqrt{5}-1}{4}d
  \]
   respectively.
The rest of the statement of this corollary also follows directly from Corollary~\ref{grp and subdiv grp spectrum}.\qed

We are going to restrict the Corollary~\ref{subdiv ped} to the case when $G$ is a regular bipartite graph with at most five eigenvalues. Moreover, we give a simpler proof of some main results proved by Kubota in~\cite{ShoPED}. 
\begin{corollary}[Theorem $3.3$, Theorem $4.1$ in~\cite{ShoPED}]
Let $G$ be a regular bipartite graph $G$ with at most five distinct eigenvalue
and moreover, one of the followings holds:
 \begin{enumerate}[label=(\alph*)]
 \item $G$ is a complete bipartite graph $K_{d,d}$, 
 \item $G$ is $C_6$, 
 \item $G$ has exactly five eigenvalue $\big\{0,\pm \frac{d}{2},\pm d\big\}$ with even $d$.
 \end{enumerate}
 \end{corollary}
\proof
A regular bipartite graph that has two or three distinct eigenvalues is a complete bipartite graph. Now consider the case when $G$ has exactly four distinct eigenvalues.  Cvetkovi\'{c} et al.~in~\cite[Page 116]{cvetkovic1980spectra} proves that a connected bipartite regular graph with four distinct eigenvalues must be the incidence graph of a symmetric $2$-$(v,d,\lambda)$ design and its spectrum is 
\[
[d]^1,\quad [\sqrt{d-\lambda}]^{v-1},\quad [-\sqrt{d-\lambda}]^{v-1},[-d]^1.
\]
The second largest eigenvalues of $A(G)$ are $\big\{\frac{\sqrt{2}}{2}d,\frac{\sqrt{3}}{2}d, \frac{1}{2}d\big\}$. Thus, one of the following three holds:
\[
\lambda=d-\frac{1}{4}d^2
,\quad
\lambda=d-\frac{3}{2}d^2,\quad
\lambda=d-\frac{1}{2}d^2
\]
For a positive integer $\lambda$, the only equation has a integer solution is $\lambda=d-\frac{1}{4}d^2$ and the solution is $d=2$. If the bipartite graph has exactly distinct four eigenvalues, then it has diameter three. Since $d=2$ here, the only feasible graph is $C_6$.

If $G$ has five eigenvalues, then $d$ must be even and the eigenvalues must be $\{ \pm d,\pm \theta,0\}$ and $\theta\in \big\{\frac{\sqrt{2}}{2}d,\frac{\sqrt{3}}{2}d, \frac{1}{2}d\big\}$.\qed

The following corollary shows that if $G$ is a regular bipartite graph and squares of its eigenvalues are rational numbers, then Grover's walk defined over $G$ is periodic with period $k$ if and only if bipartite walk defined over $G$ is periodic with period $\tau$. Moreover, by Section~\ref{1BW=2GW}, we know $k$ is even and $\tau=\frac{k}{2}$.

\begin{corollary}
Let $G$ be a regular bipartite graph. Assume that the square of each eigenvalue of $G$ is rational. Then $G$ is periodic if and only if 
$S(G)$ is periodic.

 \end{corollary}
\proof Let $\lambda_r$ be an eigenvalue of $G$. By Theorem~\ref{biregular ped}, the graph $G$ is periodic if and only if
\[
	\lambda_r^2 \in \big\{0,\frac{1}{2}d^2, \frac{3}{4}d^2,\frac{1}{4} d^2,d^2\big\}.
	\]
	By Corollary~\ref{subdiv ped}, graph $S(G)$ is periodic\qed

\section{Examples and questions}
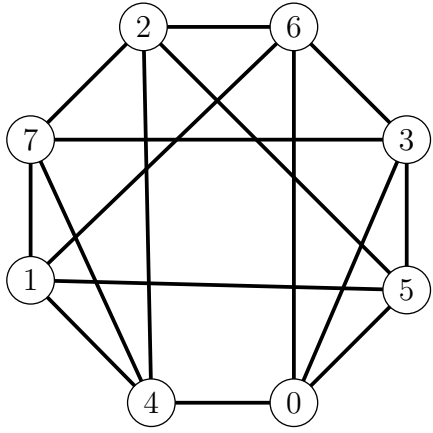
\begin{figure}[h]
\centering
\begin{tikzpicture}
\definecolor{cv0}{rgb}{0.0,0.0,0.0}
\definecolor{cfv0}{rgb}{1.0,1.0,1.0}
\definecolor{clv0}{rgb}{0.0,0.0,0.0}
\definecolor{cv1}{rgb}{0.0,0.0,0.0}
\definecolor{cfv1}{rgb}{1.0,1.0,1.0}
\definecolor{clv1}{rgb}{0.0,0.0,0.0}
\definecolor{cv2}{rgb}{0.0,0.0,0.0}
\definecolor{cfv2}{rgb}{1.0,1.0,1.0}
\definecolor{clv2}{rgb}{0.0,0.0,0.0}
\definecolor{cv3}{rgb}{0.0,0.0,0.0}
\definecolor{cfv3}{rgb}{1.0,1.0,1.0}
\definecolor{clv3}{rgb}{0.0,0.0,0.0}
\definecolor{cv4}{rgb}{0.0,0.0,0.0}
\definecolor{cfv4}{rgb}{1.0,1.0,1.0}
\definecolor{clv4}{rgb}{0.0,0.0,0.0}
\definecolor{cv5}{rgb}{0.0,0.0,0.0}
\definecolor{cfv5}{rgb}{1.0,1.0,1.0}
\definecolor{clv5}{rgb}{0.0,0.0,0.0}
\definecolor{cv6}{rgb}{0.0,0.0,0.0}
\definecolor{cfv6}{rgb}{1.0,1.0,1.0}
\definecolor{clv6}{rgb}{0.0,0.0,0.0}
\definecolor{cv7}{rgb}{0.0,0.0,0.0}
\definecolor{cfv7}{rgb}{1.0,1.0,1.0}
\definecolor{clv7}{rgb}{0.0,0.0,0.0}
\definecolor{cv0v3}{rgb}{0.0,0.0,0.0}
\definecolor{cv0v4}{rgb}{0.0,0.0,0.0}
\definecolor{cv0v5}{rgb}{0.0,0.0,0.0}
\definecolor{cv0v6}{rgb}{0.0,0.0,0.0}
\definecolor{cv1v4}{rgb}{0.0,0.0,0.0}
\definecolor{cv1v5}{rgb}{0.0,0.0,0.0}
\definecolor{cv1v6}{rgb}{0.0,0.0,0.0}
\definecolor{cv1v7}{rgb}{0.0,0.0,0.0}
\definecolor{cv2v4}{rgb}{0.0,0.0,0.0}
\definecolor{cv2v5}{rgb}{0.0,0.0,0.0}
\definecolor{cv2v6}{rgb}{0.0,0.0,0.0}
\definecolor{cv2v7}{rgb}{0.0,0.0,0.0}
\definecolor{cv3v5}{rgb}{0.0,0.0,0.0}
\definecolor{cv3v6}{rgb}{0.0,0.0,0.0}
\definecolor{cv3v7}{rgb}{0.0,0.0,0.0}
\definecolor{cv4v7}{rgb}{0.0,0.0,0.0}
\Vertex[style={minimum size=1.0cm,draw=cv0,fill=cfv0,text=clv0,shape=circle},LabelOut=false,L=\hbox{$0$},x=4cm,y=0cm]{v0}
\Vertex[style={minimum size=1.0cm,draw=cv1,fill=cfv1,text=clv1,shape=circle},LabelOut=false,L=\hbox{$1$},x=0.5cm,y=1.6292cm]{v1}
\Vertex[style={minimum size=1.0cm,draw=cv2,fill=cfv2,text=clv2,shape=circle},LabelOut=false,L=\hbox{$2$},x=2.0cm,y=5cm]{v2}
\Vertex[style={minimum size=1.0cm,draw=cv3,fill=cfv3,text=clv3,shape=circle},LabelOut=false,L=\hbox{$3$},x=5.50cm,y=3.5cm]{v3}
\Vertex[style={minimum size=1.0cm,draw=cv4,fill=cfv4,text=clv4,shape=circle},LabelOut=false,L=\hbox{$4$},x=2.1052cm,y=0.0cm]{v4}
\Vertex[style={minimum size=1.0cm,draw=cv5,fill=cfv5,text=clv5,shape=circle},LabelOut=false,L=\hbox{$5$},x=5.5cm,y=1.5cm]{v5}
\Vertex[style={minimum size=1.0cm,draw=cv6,fill=cfv6,text=clv6,shape=circle},LabelOut=false,L=\hbox{$6$},x=4cm,y=5.0cm]{v6}
\Vertex[style={minimum size=1.0cm,draw=cv7,fill=cfv7,text=clv7,shape=circle},LabelOut=false,L=\hbox{$7$},x=0.5cm,y=3.5cm]{v7}
\Edge[lw=0.05cm,style={color=cv0v3,},](v0)(v3)
\Edge[lw=0.05cm,style={color=cv0v4,},](v0)(v4)
\Edge[lw=0.05cm,style={color=cv0v5,},](v0)(v5)
\Edge[lw=0.05cm,style={color=cv0v6,},](v0)(v6)
\Edge[lw=0.05cm,style={color=cv1v4,},](v1)(v4)
\Edge[lw=0.05cm,style={color=cv1v5,},](v1)(v5)
\Edge[lw=0.05cm,style={color=cv1v6,},](v1)(v6)
\Edge[lw=0.05cm,style={color=cv1v7,},](v1)(v7)
\Edge[lw=0.05cm,style={color=cv2v4,},](v2)(v4)
\Edge[lw=0.05cm,style={color=cv2v5,},](v2)(v5)
\Edge[lw=0.05cm,style={color=cv2v6,},](v2)(v6)
\Edge[lw=0.05cm,style={color=cv2v7,},](v2)(v7)
\Edge[lw=0.05cm,style={color=cv3v5,},](v3)(v5)
\Edge[lw=0.05cm,style={color=cv3v6,},](v3)(v6)
\Edge[lw=0.05cm,style={color=cv3v7,},](v3)(v7)
\Edge[lw=0.05cm,style={color=cv4v7,},](v4)(v7)
\end{tikzpicture}
\caption{$\{-(1+\sqrt{5}),-2,0,\sqrt{5}-1,4\}$}
\label{sebi}
\end{figure}

Let $G$ denote the graph shown in Figure~\ref{sebi}. The bipartite double cover of $G$ is the Kronecker product $G\times K_2$, which is a $4$-regular bipartite graph.

The adjacency matrix of $G$ has spectrum 
\[
	\{-(1+\sqrt{5}),-2,0,\sqrt{5}-1,4\}.
	\]
The adjacency matrix of $G\times K_2$ has spectrum 
\[
\{\pm (1+\sqrt{5}),\pm 2,0,\pm(\sqrt{5}-1),\pm 4 \}.
\]

By Theorem~\ref{biregular ped}, the bipartite walk defined over $G\times K_2$ is periodic. It has period $\tau=10$.

By corollary~\ref{subdiv ped}, the bipartite walk over the subdivision graph of $G\times K_2$ is periodic. In other words, Grover's walk defined over $G\times K_2$ is periodic. Moreover, by Section~\ref{1BW=2GW}, it has period $\tau=20$.

We find the graph $G$ by accident. The only literature we have found that explicitly mentions graph $G$ is~\cite{SabiKoolen} by Cioab\u{a} et al. Proposition $4.6$ in~\cite{SabiKoolen} states that a connected $4$-regular graph with $\lambda_2>1$, then $\lambda_2\geq \sqrt{5}-1$ and the equality holds if and only $G$ is either a circulant graph $\text{Cayley}\left(\mathbb{Z}_{10},\{\pm 1,\pm 4\}\right)$ or the graph shown in Figure~\ref{sebi}. 

Note that $\text{Cayley}\left(\mathbb{Z}_{10},\{\pm 1,\pm 4\}\right)$ is a $4$-regular graph and it has eigenvalues  
\[
	\{ -(1+\sqrt{5}),0,\sqrt{5}-1,4\}.
	\]
 Let $\lambda$ be an eigenvalue of the subdivision graph of $\text{Cayley}\left(\mathbb{Z}_{10},\{\pm 1,\pm 4\}\right)$. By Corollary~\ref{subdiv ped}, the bipartite walk defined over the subdivision graph of $\text{Cayley}\left(\mathbb{Z}_{10},\{\pm 1,\pm 4\}\right)$ is periodic with period $\tau=20$.

$G\times K_2$ and its subdivision graph and the subdivision graph of $\text{Cayley}\left(\mathbb{Z}_{10},\{\pm 1,\pm 4\}\right)$ and $G\times K_2$ are only three periodic bipartite graphs we have found such that they have eigenvalues $\lambda$ with  \[ \lambda^2\in 
	\Big\{\frac{5\pm\sqrt{5}}{8}d_0d_1 ,\frac{3\pm\sqrt{5}}{8}d_0d_1\Big\}.
	\]
 We want to know how to construct periodic biregular bipartite graphs such that the squares of some of their eigenvalues that are not integers and satisfy the spectrum condition in Theorem~\ref{PED characterization}.

In Section $5$ of~\cite{ShoPED}, Kubota suggests ways to use graph products to construct regular bipartite graphs with at most five distinct eigenvalues that have periodic Grover's walks. We wonder if there are ways to construct graphs with spectrum that satisfies Theorem~\ref{PED characterization} that are not obtained from graph products.

In this paper, we study periodicity of the bipartite walk over $G$ under the assumption that squares of eigenvalues of the $G$ are algebraic integers with degree at most two. Let $X$ be a graph whose algebraic integers with degree greater than two. We want to know when the bipartite walk defined over $X$ is periodic. Author believes Theorem $3.3$ in~\cite{algtripi} will shed some light on the question.

\bibliographystyle{plain}

\end{document}